\documentclass[
   ,final           
    ,numberedheadings 
  ]
  {aipproc}

\layoutstyle{6x9}
\usepackage{hyperref}
\usepackage{amsmath,amsthm,amsfonts,amssymb}
\usepackage{array}
\numberwithin{equation}{section}

\theoremstyle{plain}  
\newtheorem{theorem}{Theorem}[section]

\newtheorem{proposition}[theorem]{Proposition}

\theoremstyle{definition}
\newtheorem{definition}[theorem]{Definition}

 \newtheorem{example}[theorem]{Example}
\newtheorem{remark}[theorem]{Remark}

\newcommand{\dbar}{\bar{\partial}}

\newcommand{\lie}{\mathfrak}

\newcommand{\R}{\mathbb{R}}
\newcommand{\Z}{\mathbb{Z}}
\newcommand{\C}{\mathbb{C}}

\newcommand{\abs}[1]{\lvert#1\rvert}

\newcommand{\x}{\times}
\newcommand{\suchthat}{\;|\;}

\newcommand{\U}{\mathrm{U}}

\newcommand{\GL}{\mathrm{GL}}
\newcommand{\SL}{\mathrm{SL}}

\newcommand{\SO}{\mathrm{SO}}
\newcommand{\Sp}{\mathrm{Sp}}

\newcommand{\Pic}{\mathrm{Pic}}
\newcommand{\End}{\mathrm{End}}
\newcommand{\rk}{\mathrm{rk}}
\newcommand{\Jac}{\mathrm{Jac}}
\newcommand{\Hom}{\mathrm{Hom}}
\newcommand{\im}{\mathrm{im}}

\newcommand{\Aut}{\mathrm{Aut}}
\newcommand{\tr}{\mathrm{tr}}
\DeclareMathOperator{\Ad}{Ad}
\DeclareMathOperator{\Id}{Id}

\newcommand{\liem}{\mathfrak{m}}

\newcommand{\liemc}{\mathfrak{m}^{\mathbb{C}}}
\newcommand{\lieh}{\mathfrak{h}}

\begin{document}

\title{Higgs bundles and the real symplectic group}



\author{Peter B. Gothen\thanks{
Member of VBAC (Vector Bundles on Algebraic Curves).
Partially supported by FCT (Portugal) through the projects
PTDC/MAT/099275/2008 and PTDC/MAT/098770/2008, and through the Centro
de Matem\'atica da Universidade do Porto.
\newline 
\emph{Date:} 06/03/2011
}}{
  address={Centro de Matem\'atica da Universidade do Porto\\
Faculdade de Ci\^encias, Universidade do Porto \\
Rua do Campo Alegre 687, 4169-007 Porto, Portugal \\}
  email: \texttt{pbgothen@fc.up.pt}
}

\begin{abstract}
We give an overview of the work of Corlette, Donaldson,
Hitchin and Simpson leading to the non-abelian Hodge theory
correspondence between representations of the fundamental group of a
surface (a \emph{surface group}) and the moduli space of Higgs
bundles. We then explain how this can be generalized to a
correspondence between character varieties for representations of
surface groups in real Lie groups $G$ and the moduli space of
$G$-Higgs bundles. Finally we survey recent joint work with Bradlow,
Garc\'\i{}a-Prada and Mundet i Riera on the moduli space of maximal
$\Sp(2n,\R)$-Higgs bundles.
\end{abstract}


\maketitle


\section{Introduction}
\label{sec:intro}

Higgs bundles are important in many areas of mathematics and
mathematical physics. For example, it was shown by Hitchin that their
moduli spaces give examples of Hyper-K\"ahler manifolds
\cite{hitchin:1987a} and that they provide important algebraically
integrable systems \cite{hitchin:1987b}. Recently they have featured
for instance in the work of Hausel--Thaddeus
\cite{hausel-thaddeus:2003} on mirror symmetry and in the work of
Kapustin--Witten \cite{kapustin-witten:2007} giving a physical
derivation of the geometric Langlands correspondence.

In this paper we start by explaining the non-abelian Hodge theory
correspondence between representations of the fundamental group of a
surface (a \emph{surface group}) and the moduli space of Higgs bundles
coming from the work of Corlette \cite{corlette:1988}, Donaldson
\cite{donaldson:1987}, Hitchin \cite{hitchin:1987a} and Simpson
\cite{simpson:1988}. 
(By the work of Corlette and Simpson the
correspondence generalizes to higher dimensional K\"ahler manifolds
but in this paper we shall deal exclusively with the case of Riemann
surfaces.) 
Next we explain how this theory can be generalized in a
systematic way to a theory of $G$-Higgs bundles for real reductive Lie
groups $G$ (giving a correspondence with surface group representations
in $G$); this is mainly based on joint work with Bradlow,
Garc\'\i{}a-Prada and Mundet i Riera.  Finally we focus on the case of
the real symplectic group $G=\Sp(2n,\R)$ and survey some recent
results on the corresponding moduli space of $\Sp(2n,\R)$-Higgs
bundles.


\section{Higgs bundles and non-abelian Hodge Theory}

\subsection{Higgs bundles}

Let $X$ be a closed Riemann surface of genus $g$ and let $K_X =
T^{1,0}X^*$ be the canonical bundle of $X$, i.e.\ its holomorphic
cotangent bundle.

\begin{definition}
  A \emph{Higgs bundle} on $X$ is a pair $(E,\phi)$, where $E \to X$
  is holomorphic vector bundle and the \emph{Higgs field} $\phi$ is a
  holomorphic $1$-form with values in $\End(E)$, i.e., $\phi \in
  H^0(X,\End(E)\otimes K_X)$.
\end{definition}

Recall that the $C^{\infty}$ isomorphism class of a complex vector
bundle $\mathcal{E}$ on the surface $X$ is given by its first Chern
class $c_1(\mathcal{E}) \in H^2(X,\Z)$ or, equivalently, by its
\emph{degree},
\begin{displaymath}
  \deg(\mathcal{E}) = \int_{X}c_1(\mathcal{E}) \in \Z.
\end{displaymath}
Recall also that $c_1(\mathcal{E})$ can be represented in de Rham
cohomology by $\frac{i}{2\pi}\tr(F(B))$, for any connection $B$ on
$\mathcal{E}$.  Thus a Higgs bundle $(E,\phi)$ has as discrete
invariants its degree $\deg(E)$ and its rank $\rk(E)$.

As a first example, a rank one Higgs bundle is a pair $(L,\phi)$, where $L
\to X$ is a line bundle and $\phi \in H^0(X,K_X)$ is a holomorphic
$1$-form (see Goldman and Xia \cite{goldman-xia:2008} for a careful
study of this case). Let $\Jac(X) = \Pic^0(X)$ be the Jacobian of $X$
which parametrizes holomorphic line bundles on $X$ of degree zero. The
tangent space to $\Jac(X)$ at any $L$ is just $H^1(X,\mathcal{O})$,
which is isomorphic to $H^0(X,K_X)^*$ by Serre duality. Hence a rank
$1$ degree zero Higgs bundle $(L,\phi)$ corresponds to a point in the
cotangent space $T^*_L\Jac(X)$ and the set of isomorphism classes of
all such $(L,\phi)$ (the \emph{moduli space}) can be identified with
the cotangent bundle to the Jacobian $T^*\Jac(X)$.
  
We now describe abelian Hodge theory correspondence for first
cohomology in a manner that points the way to the non-abelian
generalization provided by Higgs bundle theory.  Let $Y$ be a K\"ahler
manifold and let $H^{p,q}(Y)$ denote its Dolbeault cohomology groups,
defined as the $\dbar$-closed differential forms of type $(p,q)$
modulo the $\dbar$-exact ones. The abelian Hodge Theorem (see e.g.\
\cite{griffiths-harris:1978}) says that there is a decomposition of
the cohomology with complex coefficients
\begin{math}
  H^k(Y;\C) = \bigoplus_{p+q=k}H^{p,q}(Y)
\end{math}
by using harmonic representatives of cohomology classes.  
In particular, for the
Riemann surface $X$ we have
\begin{displaymath}
  H^1(X,\C) \cong H^{1,0}(X) \oplus H^{0,1}(X).
\end{displaymath}
In fact, this is the infinitesimal version of an isomorphism
\begin{equation}
  \label{eq:1}
  H^1(X,\C^*) \cong T^*\Jac(X).
\end{equation}
To see how this comes about, note that an element of
$H^1(X,\C^*)$ corresponds to a complex line bundle $\mathcal{L}$ with
a flat connection, which can be written as $d+B$ for a closed
$1$-form $B \in A^1(X,\C)$. If we fix a Hermitian metric on
$\mathcal{L}$ and use the decomposition $\C = \lie{gl}(1,\C) = i\R
\oplus \R$, we can write
\begin{displaymath}
    B = A + \psi,\qquad A \in A^1(X,i\R)\quad\text{and}\quad
      \psi\in A^1(X,\R).
\end{displaymath}
Then the flatness condition $dB=0$ is equivalent to saying that $A$
and $\psi$ are both closed $1$-forms. Thus $A$ is a flat unitary
connection on $\mathcal{L}$. Moreover, the Hodge Theorem says that the
cohomology class of $\psi$ has a harmonic representative $\psi+df$
for some smooth $f$, meaning that
\begin{displaymath}
  d^*(\psi+df)=0.
\end{displaymath}
Note that $\psi+df$ is obtained by applying to $\psi$ the gauge
transformation $\chi = \exp(f)$ of $\mathcal{L}$. Any gauge
transformation preserves the flatness condition $dB=0$ and hence the
pair $(A,\psi)$ obtained applying the above procedure to the gauge
transformed $\chi\cdot B$ satisfies the equations
\begin{align}
    dA &=0, \label{eq:2} \\
    d\psi&=0, \label{eq:3} \\
    d^*\psi&=0. \label{eq:4}
\end{align}
Now the connection $A$ defines a holomorphic line bundle $L_A\to X$ by
taking the holomorphic structure on $\mathcal{L}$ given by the
associated $\dbar$-operator $\dbar_A$, which is just the $(0,1)$-part
of the covariant derivative $d_A$. We also define 
\begin{displaymath}
\phi =
\psi^{1,0} \in \Omega^{1,0}(X).
\end{displaymath}
The important point to note is now that the conditions (\ref{eq:3}) and
(\ref{eq:4}) together are equivalent to the holomorphicity condition
$\dbar_A\phi = 0$ and hence (\ref{eq:2})--(\ref{eq:4}) are equivalent
to the pair of equations 
\begin{equation}
  \label{eq:5}
  \begin{aligned}
    dA &=0, \\
    \dbar_A\phi&=0.
  \end{aligned}
\end{equation}
Thus
\begin{math}
  \phi \in H^0(X,K_X)
\end{math}
is holomorphic and $(L_A,\phi)$ is a rank $1$ Higgs bundle. This, as
we have seen, corresponds to a point in $T^*\Jac(X)$ (the fact that
$\deg(L_A)=0$ follows from $F(A) = dA = 0$).

An analogous argument shows how to recover a flat line bundle from a
rank $1$ Higgs bundle $(L,\phi)$, completing the proof of
(\ref{eq:1}).

\begin{remark}
  \label{rem:abelian-character-variety}
  It follows that the moduli space $T^*\Jac(X)$ of degree zero Higgs
  line bundles can be identified with the variety of complex
  characters of $\pi_1(X)$ via the isomorphism $H^1(X,\C^*) \cong
  \Hom(\pi_1(X),\C^*)$.  However the natural complex structures on
  these spaces do not correspond under this identification.
\end{remark}

\subsection{Non-abelian Hodge Theory}

In this section we explain how the abelian Hodge Theorem for first
cohomology generalizes to higher rank.  Let $\mathcal{E}\to X$ be a
fixed $C^\infty$ complex vector bundle of rank $n$ and degree $d$.
Let $B$ be a connection in $\mathcal{E}$ with constant central
curvature (such a connection is also called \emph{projectively flat}):
\begin{equation}
  \label{eq:12}
  F(B) = -i\mu\Id\omega,
\end{equation}
where $\omega$ is the K\"ahler class of $X$, normalized so that
$\int_X \omega = 2\pi$. Taking the trace and integrating in this
formula, we see that the constant $\mu$ is given by
\begin{displaymath}
  \mu = \frac{\deg(\mathcal{E})}{\rk(\mathcal{E})} = \frac{n}{d};
\end{displaymath}
the ratio $\mu(E) = \frac{\deg(\mathcal{E})}{\rk(\mathcal{E})}$ is
called the \emph{slope} of $\mathcal{E}$.
In particular, if $\deg(\mathcal{E}) = 0$, the connection $B$ has zero
curvature, i.e., it is flat.

Let $h$ be a
Hermitian metric in $\mathcal{E}$. We shall write
$\mathcal{E}_{h}$ for the bundle $\mathcal{E}$ equipped with Hermitian
metric $h$. Splitting the associated Lie
algebra valued $1$-form in its skew-Hermitian and Hermitian parts, we
can write
\begin{displaymath}
  B = A + \psi,
\end{displaymath}
where $A$ is a unitary connection on $\mathcal{E}$ and $\psi$ takes
values in the bundle $\End^{\mathrm{herm}}(\mathcal{E}_h)$ of Hermitian
endomorphisms of $\mathcal{E}$, in other words, $\psi$ descends to
a $1$-form
\begin{displaymath}
  \psi \in \Omega^1(X,\End^{\mathrm{herm}}(\mathcal{E}_h)).
\end{displaymath}
By taking the skew-Hermitian and Hermitian parts of (\ref{eq:12}), we
can express this condition in terms of $(A,\psi)$ as follows:
\begin{align}
  F(A) + \tfrac{1}{2}[\psi,\psi] &=-i\mu\Id\omega, \label{eq:6} \\
  d_A\psi &=0.
\end{align}
One may ask whether there is a preferred choice of Hermitian metric
in the projectively flat bundle $\mathcal{E}$. Since the space of metrics in
$\C^n$ is $\GL(n,\C) / \U(n)$, a Hermitian metric $h$ may be viewed as a
$\pi_1(X)$-equivariant map
\begin{displaymath}
  \tilde h\colon \tilde{X} \to \GL(n,\C) / \U(n).
\end{displaymath}
Because the symmetric space $\GL(n,\C) / \U(n)$ is a (negatively
curved) Riemannian manifold and we have a conformal class of metrics
on $X$ defined by the complex structure, it makes sense to ask for
$\tilde{h}$ to be an equivariant harmonic\footnote{Recall that
  harmonicity only depends on the conformal class of the source
  manifold when this is $2$-dimensional.} map. 
The derivative of
$\tilde h$ is a section
\begin{displaymath}
  d\tilde h \in \Omega^1(\tilde{X}, \tilde{h}^*T(\GL(n,\C) / \U(n))).
\end{displaymath}
The harmonic map equation (cf. \cite{eells-sampson:1964}) for $\tilde h$ is 
\begin{equation}
  \label{eq:8}
   d_{\nabla}^*(d\tilde{h}) = 0,
\end{equation}
where $d_{\nabla}$ is the pull-back by $h$ of the Levi-Civita
connection on $\GL(n,\C) / \U(n))$ and $d_{\nabla}^*$ is its adjoint
(constructed again using the conformal class of metrics on $X$ and the
Riemannian metric on $\GL(n,\C) / \U(n))$).  Moreover, it can be shown
that (\ref{eq:8}) is equivalent to the equation
\begin{math}
  d_A^*\psi = 0
\end{math}
so that the pair $(A,\psi)$ obtained from $B$ via a harmonic
metric satisfies the equations 
\begin{align}
  F(A) + \tfrac{1}{2}[\psi,\psi] &=-i\mu\Id\omega, \label{eq:9}\\
  d_A\psi &=0, \label{eq:10}\\
  d_A^*\psi &=0. \label{eq:7}
\end{align}
    
Before we can state the main existence result for harmonic metrics we
need to introduce the notion of a reductive  connection.
A connection $B$ on $\mathcal{E}$ corresponds to a covariant
derivative $d_B\colon \Omega^0(X,\mathcal{E}) \to
\Omega^1(X,\mathcal{E})$. We say that a subbundle $\mathcal{F} \subset
\mathcal{E}$ is preserved by $B$ if it satisfies
$d_B(\Omega^0(X,\mathcal{F})) \subset \Omega^1(X,\mathcal{F})$.

\begin{definition}
  A connection $B$ on $\mathcal{E}$ is \emph{reductive} if for any
  subbundle $\mathcal{F} \subset \mathcal{E}$ preserved by $B$, there
  is a subbundle $\mathcal{F}'$ which is also preserved by $B$ and
  such that $\mathcal{E} = \mathcal{F} \oplus \mathcal{F}'$.
\end{definition}

It is not hard to show that if a bundle $\mathcal{E}$ with a constant
central curvature connection $B$ admits a harmonic metric, then $B$ is
reductive. The following theorem says that the converse holds. It is
due to Donaldson~\cite{donaldson:1987} (in the case of rank $2$
bundles on a Riemann surface) and Corlette \cite{corlette:1988} (for
base manifolds of higher dimension and more general structure groups).

\begin{theorem}
  Let $B$ be a connection of constant central curvature in
  $\mathcal{E}$. Then there exists a unique harmonic metric in
  $\mathcal{E}$.
\end{theorem}

In order to get a global statement, we introduce the \emph{moduli
  space of reductive connections of constant central curvature} on a
 vector bundle $\mathcal{E}$ of degree $d$ and rank $n$:
\begin{displaymath}
  M^{\mathrm{dR}}_d(X,\GL(n,\C)) =
  \left\{ B \suchthat \text{$F(B)=-i\mu\Id\omega$ and $B$ is reductive}\right\} 
    / \mathcal{G},
\end{displaymath}
where $\mathcal{G}=\Aut(\mathcal{E})$ is the complex gauge group of
$\mathcal{E}$. 

We remark that connections of constant central curvature are related to
representations of the fundamental group of $X$ as follows (see Atiyah--Bott
  \cite{atiyah-bott:1982}).
There is a central extension
\begin{displaymath}
  0 \to \R \to \Gamma_{\R} \to \pi_1 X \to 1 
\end{displaymath}
defined by $\Gamma_{\R} = \Gamma\otimes_{\Z}\R$, where the universal
central extension $\Gamma$ is defined by
\begin{displaymath}
  \Gamma = \langle a_1,b_1,\dots,a_g,b_g,J \suchthat
  \text{$J$ is central and $\prod[a_i,b_i] = 1$} \rangle.
\end{displaymath}
Define the \emph{character variety} for representations of $\Gamma_{\R}$
in $\GL(n,\C)$ by
\begin{multline*}
  M^{B}_d(X,\GL(n,\C)) \\
    = \{\rho\colon\Gamma_{\R}\to\GL(n,\C) \suchthat
      \text{$\rho(J)
      = \exp(2\pi id/n)\Id$ and $\rho$ is semisimple}\} / \GL(n,\C),
\end{multline*}
where $\GL(n,\C)$ acts by overall conjugation. Note that
$M^{B}_0(X,\GL(n,\C))$ can be naturally identified with the character
variety for representations of $\pi_1 X$ in $\GL(n,\C)$.
It is now a standard fact that there is an (analytic) isomorphism 
\begin{equation}
  \label{eq:13}
  M^{\mathrm{dR}}_d(X,\GL(n,\C)) \cong M^{B}_d(X,\GL(n,\C)) 
\end{equation}
obtained by taking a connection $B$ to its holonomy representation.


Now fix a reference metric $h_0$ in $\mathcal{E}$.  For a reductive
connection $B$ with constant central curvature, let $h$ be the
harmonic metric given by the Theorem, so that $(\tilde
A,\tilde\psi)$ obtained from $B$ using $h$ satisfies
(\ref{eq:7}). Let $g\in \mathcal{G}$ be an isometry between unitary
bundles $g\colon\mathcal{E}_{h_0} \to \mathcal{E}_{h}$.  Then the
connection $g^*B$ also has constant central curvature and the 
pair $(A,\psi)=g^*(\tilde A,\tilde\psi)$ obtained from $g^*B$
using the metric $h_0$, will also solve the equation
(\ref{eq:7}). In other words, the metric $h_0$ is harmonic for
$g^*B$. Any two choices of $g$ differ by a unitary gauge
transformation of $(\mathcal{E},h_0)$. Hence, letting $\mathcal{U}$
denote the unitary gauge group, the theorem can be reformulated as
follows:

\begin{theorem}
  \label{thm:DC}
  There is a bijective correspondence
  \begin{displaymath}
    M^{\mathrm{dR}}_d(X,\GL(n,\C)) \cong \{(A,\psi) \suchthat
      \text{satisfying (\ref{eq:9})--(\ref{eq:7})}\}/ \mathcal{U}
    \end{displaymath}
  where $A$ is a unitary connection on $(\mathcal{E},h_0)$ and
  $\psi\in \End^{\mathrm{herm}}(\mathcal{E}_{h_0})$. 
\end{theorem}

In order to obtain a Higgs bundle from a solution $(A,\psi)$
to (\ref{eq:9})--(\ref{eq:7}), we decompose the covariant derivative
associated to $A$ as
\begin{displaymath}
  d_A = \partial_A + \dbar_A,
\end{displaymath}
and denote by $E_A$ the holomorphic
vector bundle defined by $\dbar_A$. Similarly, one can write
\begin{displaymath}
  \psi = \phi+\phi^*
\end{displaymath}
for a unique $\phi \in \Omega^{1,0}(X,\End(\mathcal{E}))$ (here
$\phi^*$ denotes the $(0,1)$-form obtained from $\phi$ by combining
the adjoint with respect to $h_0$ with conjugation on the form component).
Then one easily checks that (\ref{eq:9})--(\ref{eq:7}) are
equivalent to \emph{Hitchin's equations}
\begin{equation}
\label{eq:11}
\begin{aligned}
  F(A) + [\phi,\phi^*]& = -i\mu\Id\omega, \\
  \dbar_A\phi &=0.
\end{aligned}
\end{equation}
Note that the second equation says that $\phi$ is a holomorphic
$1$-form so $(E_A,\phi)$ is a Higgs bundle. Now recall that a
holomorphic vector bundle which admits a flat unitary connection is
the direct sum of stable degree zero vector bundles.
To make sense of
the analogous statement for Higgs bundles, we define the
following stability
notions.
\begin{definition}
  A Higgs bundle $(E,\phi)$ is
    \begin{itemize}
    \item \emph{semistable} if $\mu(F) \leq \mu(E)$ for all
      holomorphic subbundles
      $F \subset E$ such that $\phi(F) \subset F \otimes K_X$;

    \item \emph{stable} if $\mu(F) < \mu(E)$ for all non-zero proper
      holomorphic subbundles $F \subsetneq E$ such that $\phi(F)
      \subset F \otimes K_X$;

    \item \emph{polystable} if $(E,\phi) = (E_1,\phi_1) \oplus \dots
      \oplus (E_r,\phi_r)$, where each $(E_i,\phi_i)$ is stable with
      $\mu(E_i) = \mu(E)$.
    \end{itemize}
\end{definition}
A subbundle $F \subset E$ such that $\phi(F) \subset F \otimes K_X$
is said to be a \emph{$\Phi$-invariant subbundle}. Note that
semistability is a weaker condition than polystability, which in turn
is weaker than stability.

It is easy to check that 
a Higgs bundle obtained from a solution to
Hitchin's equations (\ref{eq:11}) is polystable.  The
converse is given by the following theorem, which gives a
Hitchin--Kobayashi correspondence for Higgs bundles. It is due to
Hitchin \cite{hitchin:1987a} (for Higgs bundles on Riemann surfaces)
and Simpson \cite{simpson:1988} (for Higgs bundles over higher
dimensional manifolds).

\begin{theorem}
  \label{thm:HS}
  If $(E,\phi)$ is polystable then there exists a unique Hermitian
  metric in $E$ such that $(A,\phi)$ satisfies Hitchin's equations
  (\ref{eq:11}), where $A$ is the unique unitary connection compatible
  with the holomorphic structure (i.e.\ the Chern connection).
\end{theorem}

In order to get the corresponding global statement, we introduce 
the \emph{moduli space} 
\begin{displaymath}
  M^{\mathrm{Dol}}_d(X,\GL(n,\C))
\end{displaymath}
of rank $n$, degree $d$ polystable Higgs bundles. As a set, this is
the set of isomorphism classes of polystable Higgs bundles. It can be
given the structure of a complex (algebraic) variety using standard
gauge theory methods (Hitchin \cite{hitchin:1987a}) or using Geometric
Invariant Theory\footnote{From this point of view it is better to
  consider $M^{\mathrm{Dol}}_d(X,\GL(n,\C))$ as the space of
  $S$-equivalence classes of semistable Higgs bundles.} (Nitsure
\cite{nitsure:1991}).  Then Theorem~\ref{thm:HS} implies that
$M^{\mathrm{Dol}}_d(X,\GL(n,\C))$ is in bijective correspondence with
the space of unitary gauge equivalence classes of solutions to
Hitchin's equations (\ref{eq:11}). Putting this together with
Theorem~\ref{thm:DC} and the identification (\ref{eq:13}), we finally
obtain the non-abelian Hodge Theorem.

\begin{theorem}
    There is a homeomorphism
  \begin{math}
    M^{\mathrm{B}}_d(X,\GL(n,\C)) \cong M^{\mathrm{Dol}}_d(X,\GL(n,\C)).
  \end{math}
\end{theorem}

\section{$G$-Higgs bundles}

We have seen that Higgs bundles correspond to representations $\rho\colon\pi_1
X\to\GL(n,\C)$ The use of Higgs bundle methods for studying
representations
\begin{displaymath}
  \rho\colon \pi_1X \to G
\end{displaymath}
for more general Lie groups $G$ was pioneered by Hitchin
\cite{hitchin:1987a,hitchin:1992} and also, using Tannakian
considerations, by Simpson \cite{simpson:1992}.
Subsequently a theory of $G$-Higgs bundles, appropriate for studying
representations of $\pi_1 X$ in real a reductive Lie group $G$, has
been developed in a systematic way. In this section we briefly outline
this theory. For more details the reader may consult, for example,
\cite{bradlow-garcia-prada-gothen:2003,garcia-prada-gothen-mundet:2008,garcia-gothen-mundet:2009b,bradlow-garcia-prada-gothen:hss-higgs}.

Let $G$ be a real reductive Lie group in the sense of Knapp
\cite[p.~384]{knapp:1996}. In particular this means that we are given
a maximal compact subgroup $H \subset G$. 
Also, there is a \emph{Cartan decomposition}
\begin{displaymath}
\lie{g} = \lieh \oplus \liem,
\end{displaymath}
where $\lieh$ is the Lie algebra of $H$. Moreover, restriction to $H$
of the adjoint action of $G$ on its Lie algebra gives a representation
\begin{align*}
  \iota\colon H &\to \GL(\lie{m}) \\
  g &\mapsto \bigl(x\mapsto \Ad(g)(x)\bigr).
\end{align*}
This representation is called the \emph{isotropy representation}. We
shall denote by the same symbol its complexification $\iota\colon
H^{\C}\to \GL(\liemc)$ defined on the complexification $H^{\C}$ of
$H$. If $E$ is a principal $H^{\C}$-bundle, we thus have an associated
bundle 
$$E(\lie{m}^{\C}) = E \x_{\iota} \lie{m}^{\C}$$
with fibres $\liemc$. Note that if $G$ is itself a complex group, then
$\lie{g} = \lie{h} \oplus i \lie{h}$ so that $\lie{m}^{\C}=
\lie{h}^{\C} = \lie{g}$, and hence $E(\lie{m}^{\C}) = E(\lie{g})$, the
adjoint bundle of the $G$-bundle $E$.

\begin{definition}
  A \emph{$G$-Higgs bundle} is a pair $(E,\phi)$, where
 $E \to X$ is a holomorphic principal $H^{\C}$-bundle,
 $\phi \in H^0(X, E(\lie{m}^{\C}) \otimes K_X)$.
\end{definition}

\begin{example}
  Let $G = \GL(n,\C)$. Then a $G$-Higgs bundle gives rise to a Higgs
  bundle $(E\x_{\GL(n,\C)}\C^n,\phi)$ as previously defined.
\end{example}

\begin{example}
  \label{ex:sp2nr}
  Let $G = \Sp(2n,\R)$. In this case: $H = \U(n)$ and $\liemc =
  S^2\C^n \oplus S^2(\C^n)^*$, where $\C^n$ denotes the fundamental
  representation of $H^{\C} = \GL(n,\C)$.  Hence a $\Sp(2n,\R)$-Higgs
  bundle is equivalent to a triple $(V,\beta,\gamma)$, where $V$ is a
  rank $n$ vector bundle and
    \begin{displaymath}
      \beta \in H^0(X, S^2V\otimes K_X),\qquad \gamma\in
      H^0(X,S^2V^*\otimes K_X).
    \end{displaymath}
    The usual Higgs bundle given by the inclusion $\Sp(2n,\R) \subset
    \Sp(2n,\C) \subset \GL(2n,\C)$ is:
    \begin{displaymath}
      (E=V \oplus V^*,\phi=
      \begin{pmatrix}
        0 & \beta \\
        \gamma & 0
      \end{pmatrix}).
    \end{displaymath}
\end{example}


\begin{example}
  \label{ex:glnr}
  A $\GL(n,\R)$-Higgs bundle is given by
  \begin{math}
    ((W,Q),\eta)
  \end{math}, where $(W,Q)$ is an orthogonal bundle and $\eta\colon
  W \to W \otimes K_X$ is symmetric with respect to $Q$.
\end{example}

Stability of $G$-Higgs bundles is in general a complicated notion, and
we shall not state it here, since we shall have no explicit need for
it. The interested reader is referred to
\cite{garcia-gothen-mundet:2009b}. But it is worth remarking that in
all cases of interest to us here, the complexification $G^{\C} \subset
\GL(n,\C) $ is a linear group and semi- and polystability of
$(E,\phi)$ is equivalent to semi- and polystability of the induced
rank $n$ usual Higgs bundle, respectively (cf.\
\cite{bradlow-garcia-prada-gothen:2003}). On the other hand, the
stability conditions are, in general, different.

When $G$ is connected, the topological classification of $G$-bundles
is given by a characteristic class in $H^2(X,\pi_1G) \cong \pi_1G
\cong \pi_1 H$. For a fixed topological class $d \in \pi_1 H$ we
can introduce the analogues of the moduli spaces defined above. Thus
$M^{\mathrm{B}}_d(X,G)$ denotes the character variety for
representations of a suitable central extension of $\pi_1X$ and
$M^{\mathrm{Dol}}_d(X,G)$ denotes the moduli space of polystable
$G$-Higgs bundles. In order to construct the latter space, the general
theory of moduli of decorated bundles of Schmitt \cite{schmitt:2008}
is required.

We have the following generalization of the non-abelian Hodge Theorem,
proved via an intermediate moduli space of solutions to an appropriate
version of Hitchin's equations (\ref{eq:11}) for $G$-Higgs
bundles. While Theorem~\ref{thm:DC} essentially applies unchanged in
this situation, the generalization of Theorem~\ref{thm:HS} to
principal pairs of
\cite{mundet:2000,bradlow-garcia-prada-mundet:2003,garcia-gothen-mundet:2009b}
is required for proving this result.

\begin{theorem}
  There is a homeomorphism
  \begin{math}
    M^{\mathrm{B}}_d(X,G) \cong M^{\mathrm{Dol}}_d(X,G).
  \end{math}
\end{theorem}

\section{The real symplectic group}

In this section we mainly focus on the moduli space of
$\Sp(2n,\R)$-Higgs bundles, describing some properties of their moduli
spaces. Details of these results can be found in
\cite{gothen:2001,garcia-prada-mundet:2004,garcia-gothen-mundet:2009a,bradlow-garcia-gothen-mundet:2009}.

\subsection{Stability and the Milnor--Wood inequality}

Let $(V,\beta,\gamma)$ be a $\Sp(2n,\R)$-Higgs bundle as in
Example~\ref{ex:sp2nr}. The topological invariant of
$(V,\beta,\gamma)$ is the degree $d=\deg(V)$. If $(V,\beta,\gamma)$ is
polystable and $\rho\colon\pi_1 X \to \Sp(2,\R)$ is the 
corresponding representation, it can be seen that $d$ is the
so-called \emph{Toledo
  invariant} of $\rho$, usually denoted by $\tau(\rho)$. 
The Toledo invariant is bounded by the inequality
\begin{equation}
  \label{eq:16}
  \abs{\tau(\rho)} \leq n(g-1),
\end{equation}
usually known as the \emph{Milnor--Wood inequality}. This inequality is due to
Milnor \cite{milnor:1957} in the case $n=1$ and
to Dupont \cite{dupont:1978} and Turaev \cite{turaev:1984} in the
general case, the latter giving 
the sharp bound. This inequality can be proved easily using
$\Sp(2n,\R)$-Higgs bundles as follows.

For definiteness, assume that $\deg(V)> 0$. 
Define subbundles $N \subset V$ and $I \subset V^*$ using the
subsheaves $\ker(\gamma)$ and $\im(\gamma)\otimes K_X^{-1}$. Then
    \begin{displaymath}
      N \oplus 0 \subset V \oplus V^*\qquad\text{and}\qquad
      V \oplus I \subset V \oplus V^*
    \end{displaymath}
    are $\Phi$-invariant subbundles of the $\GL(2n,\R)$-Higgs bundle
    $(E= V \oplus V^*, \Phi = \left(
       \begin{smallmatrix}
         \beta & 0 \\
         0 & \gamma
       \end{smallmatrix}
     \right))$.
If $(V,\beta,\gamma)$ is semistable, then so is $(E,\Phi)$. It follows
that
\begin{align}
  \label{eq:14}
  \deg(N) &\leq 0, \\
  \deg(V) + \deg(I) &\leq 0.
\end{align}
Moreover, $\gamma\neq 0$ since otherwise $V \subset E$ would be
$\Phi$-invariant and thus violate semistability. Hence it induces a
non-zero section
\begin{equation}
  \label{eq:15}
  \bar{\gamma} \in H^0(X,\det(V/N)^*\otimes\det(I)\otimes K^{\rk(\gamma)}),
\end{equation}
so that this line bundle must have positive degree. Together with
(\ref{eq:14}) this implies that
\begin{displaymath}
  \deg(V) \leq \rk(\gamma)(g-1),
\end{displaymath}
thus demonstrating that $\deg(V) \leq n(g-1)$, which proves the
Milnor--Wood inequality (\ref{eq:16}) for $\deg(V)=\tau(\rho)>0$. An
analogous argument using $\beta$ proves the case $\deg(V)<0$. 

Note that when the equality $\deg(V) = \rk(\gamma)(g-1)$ holds, the
line bundle in (\ref{eq:15}) is forced to have degree zero. Thus the
map $\bar{\gamma}$ in (\ref{eq:15}) is an isomorphism. In particular, we
have the following important consequence.

\begin{proposition}
  \label{prop:max-gamma-iso}
  Let $(V,\beta,\gamma)$ be a semistable $\Sp(2n,\R)$-Higgs bundle
  such that $\deg(V)=n(g-1)$. Then $\gamma$ induces an isomorphism
  \begin{math}
    \gamma\colon V \xrightarrow{\cong} V \otimes K.
  \end{math}
  If $\deg(V) = -n(g-1)$ the analogous statement holds for $\beta$.
\end{proposition}
  
\begin{definition}
  An $\Sp(2n,\R)$-Higgs bundle $(V,\beta,\gamma)$ is said
  to be \emph{maximal} if $\abs{\deg(V)}=n(g-1)$. 
  Similarly, a representation $\rho\colon\pi_1 X \to \Sp(2n,\R)$ is
  \emph{maximal} if $\abs{\tau(\rho)}=n(g-1)$.
\end{definition}

The geometric importance of maximal representations is underlined by
the following theorem, due to Goldman \cite{goldman:1980}.
\begin{theorem}
    $\rho\colon \pi_1 X \to \Sp(2,\R)$ is maximal if and only if it is
    Fuchsian.
\end{theorem}
Recall that $\Sp(2,\R)=\SL(2,\R)$ and that a representation
$\rho\colon\pi_1 X\to\SL(2,\R)$ is called Fuchsian if it is discrete
and faithful. More generally, it was shown by Burger--Iozzi--Wienhard
\cite{burger-iozzi-wienhard:2003,burger-iozzi-wienhard:2010} that
maximal representations are discrete, faithful and reductive.

\subsection{The moduli space of maximal $\Sp(2n,\R)$-Higgs bundles}

For brevity, denote by $M_{\mathrm{max}}(n)$ the moduli
space $M_{n(g-1)}^{\mathrm{Dol}}(X,\Sp(2n,\R))$ of maximal $\Sp(2n,\R)$-Higgs bundles.  In this
section we shall see how Proposition~\ref{prop:max-gamma-iso} leads to
the existence of new invariants of maximal $\Sp(2n,\R)$-Higgs
bundles. In particular, this implies that
$M_{\mathrm{max}}(n)$ is disconnected and we
shall also give a complete count of its connected components.

Choose a square root $K_X^{1/2}$ of $K_X$. Using the isomorphism
$\gamma\colon V\xrightarrow{\cong}V^*\otimes K_X$ and the fact that
$\gamma$ is symmetric, we define an orthogonal holomorphic bundle
$(W,Q)$ on $X$ as follows:
\begin{displaymath}
  W = V \otimes K_X^{-1/2}\qquad\text{and}\qquad
  Q = \gamma \otimes 1_{K^{-1/2}}\colon W \xrightarrow{\cong} W^*.
\end{displaymath}
The Stiefel--Whitney classes of $(W,Q)$ define new invariants of
the $\Sp(2n,\R)$-Higgs bundle $(V,\beta,\gamma)$:
\begin{displaymath}
 w_1(V,\beta,\gamma) \in H^1(X,\Z/2)\qquad\text{and}\qquad
 w_2(V,\beta,\gamma) \in H^2(X,\Z/2).
\end{displaymath}

\begin{remark}
  In the case $n=1$, the invariant $w_2(V,\beta,\gamma)$ always
  vanishes for obvious reasons.  The case $n=2$ is also special: when
  $w_1(V,\beta,\gamma) = 0$, there is a lift of $w_2(V,\beta,\gamma)$
  to an invariant $c(V,\beta,\gamma) \in H^2(X,\Z) \cong \Z$, coming
  from a reduction of structure group to $\SO(2) \subset
  \mathrm{O}(2)$. This invariant satisfies
  $\abs{c(V,\beta,\gamma)}\leq 2g-2$.
\end{remark}

\begin{remark}
  \label{rem:hss}
  We can define $\eta = (\beta\otimes 1) \circ (\gamma\otimes
  1)\colon W \to W \otimes K_X^2$.  Then $((W,Q),\eta)$ is a
  \emph{twisted} $\GL(n,\R)$-Higgs bundle (cf.\
  Example~\ref{ex:glnr}), the difference to a usual $\GL(n,\R)$-Higgs
  bundle being the twisting by $K_X^2$ rather than $K_X$.

  This is an instance of a general phenomenon occurring for maximal
  $G$-Higgs bundles when $G$ is isogenous to the isometry group of
  a Hermitian symmetric space of non-compact type. We refer to
  \cite{bradlow-garcia-prada-gothen:hss-higgs} for this theory.
\end{remark}

The connected components of
$M_{\max}(1)=M^{\mathrm{B}}_{g-1}(X,\SL(2,\R))$ were determined by
Goldman \cite{goldman:1988} 
working directly with representations of $\pi_1 X$ (see Hitchin
\cite{hitchin:1987a} for a proof based on Higgs bundle theory):
\begin{theorem}
  The connected components of $M_{\max}(1)$ are the $2^{2g}$ subspaces
  \begin{math}
    M_{w_1} \subset M_{\max}(1)
  \end{math}
  of Higgs bundles having invariant $w_1 \in H^1(X, \Z/2)$.
\end{theorem}

\begin{remark}
  \label{rem:hitchin-component}
  These components are all homeomorphic to the Teichm\"uller space of
  $X$, as also follows from results of Goldman.  In particular each
  component ${M}_{w_1}$ is homeomorphic to $\R^{6g-6}$. This
  can bee easily seen from the Higgs bundle point of view (see
  \cite{hitchin:1987a}) by noting that a maximal $\Sp(2,\R)$-Higgs
  bundle is isomorphic to one of the form $(L,\beta,\gamma)$, where
  $L^2 = K_X$ and $\gamma=1$. The choice of $L$ given by $w_1$ and
  thus the choice of $\beta\in H^0(X,K_X^2)$ gives a an identification
  ${M}_{w_1} \cong H^0(X,K_X^2)$.
\end{remark}

Generalizing the parametrization by quadratic differentials of
${M}_{w_1}$ of the preceding remark, Hitchin
\cite{hitchin:1992} showed the
existence of special connected \emph{Hitchin components}
${M}^H \subset {M}^{\mathrm{Dol}}(X,G)$
whenever $G$ is a split real form of a simple complex group (the classical
examples are $G=\SL(n,\R)$, $\Sp(2n,\R)$, $\SO(n,n)$, $\SO(n+1,n)$).
The Hitchin components are vector spaces of the form
\begin{math}
  {M}^H \cong \bigoplus H^0(X,K_X^{m_i+1}).
\end{math}
In the case $G=\Sp(2n,\R)$ Hitchin components are maximal and there
are $2^{2g}$ such components
\begin{displaymath}
M^H_L \subset M_{\max}(n)
\end{displaymath}
indexed by square roots
$L$ of $K_X$, just as in the case $n=1$.


Denote by $M_{w_1,w_2} \subset M_{\max}(n)$ the subspace of
\textbf{non}-Hitchin $\Sp(2n,\R)$-Higgs bundles with invariants
$w_i\in H^i(X,\Z/2)$ for $i=1,2$. In the case $n=2$, we additionally
write $M_{0,c}$ for the subspace of non-Hitchin $\Sp(4,\R)$-Higgs
bundles with invariants $w_1=0$ and $c\in H^2(X,\Z) \cong \Z$.

The connected components of ${M}_{\max}(n)$ were determined in
\cite{gothen:2001} for $n=2$ and in \cite{garcia-gothen-mundet:2009a}
for $n \geq 3$. We refer the interested reader to these papers for the
proof of the following theorem. For more information on maximal
$\Sp(2n,\R)$-Higgs bundles and the corresponding representations see
for instance
\cite{bradlow-garcia-gothen-mundet:2009} and \cite{guichard-wienhard:2010}.

\begin{theorem}
  For $n=2$, the decomposition in connected components of ${M}_{\max}$
  is
  \begin{displaymath}
    {M}_{\max} = 
    \bigcup_{w_1 \neq 0, w_2} {M}_{w_1,w_2}
    \cup\bigcup_{0\leq c <2g-2} {M}_{0,c}
    \cup\bigcup_{L^2=K_X} {M}^{H}_{L}.
  \end{displaymath}
  For $n\geq 3$, the decomposition in connected components of
  ${M}^{\max}$ is
    \begin{displaymath}
      {M}_{\max} = 
      \bigcup_{w_1,w_2} {M}_{w_1,w_2}
      \cup\bigcup_{L^2=K_X} {M}^{H}_{L}.
  \end{displaymath} 
\end{theorem}

Finally, we mention that the maximal connected components has been
carried our in many cases for many non-compact groups $G$ of Hermitian
type; see \cite{bradlow-garcia-prada-gothen:hss-higgs} for a survey of
such results. On the other hand, the determination of non-maximal
components is in general a difficult problem, which in the case of
$G=\Sp(2n,\R)$ has only been carried out for $n=1$ (by Hitchin
\cite{hitchin:1987a} and Goldman \cite{goldman:1988}) and for $n=2$
\cite{garcia-prada-mundet:2004}.


\begin{theacknowledgments}
  This paper is based on joint work with S. Bradlow,
  O. Garc\'\i{}a-Prada and I. Mundet i Riera, to whom I am greatly
  indebted. I also wish to thank O. Guichard for useful comments on an
  earlier version of the paper. 
\end{theacknowledgments}



\end{document}